# VLSI Layouts and DNA Physical Mappings

Michael J. Dinneen

*Computer Research and Applications*

*Los Alamos National Laboratory*

*Los Alamos, N.M. 87545*

## 1  Introduction

In this short note, we show that an important problem in computational biology is equivalent to a colored version of a well-known graph layout problem. In order to map the human genome, biologists use graph theory, particularly interval graphs, to model the overlaps of DNA clones (cut up segments of a genome) [Mir94]. For engineers, Very-Large-Scale-Integrated (VLSI) circuits must be laid out in order to minimize physical and cost constraints. The vertex separation (see below) of a graph layout is one such measurement of how good a layout is.

The NP-complete combinatorial problem of *Intervalizing Colored Graphs* (ICG) first defined in [FHW93] (and independently given in [GKS93] as the *Graph Interval Sandwich* problem) is intended to be a limited, first-step model for finding DNA physical mappings. For this model, it is assumed that the biologist knows some of the overlaps — for instance, overlaps specified by some probability threshold based on the physical data. The question asked by the ICG problem is whether other edges can be properly added to differently colored vertices to form a colored interval graph.

Finding the *Vertex Separation* (VS) of a graph is related to many diverse problems in computer science besides its importance to VLSI layouts. Lengauer showed that progressive black/white pebble game (important to compiler theory) and vertex separation are polynomially reducible to each other [Len81]. Node search number, a variant of search number [Par76], was shown equivalent to the vertex separation plus one by Kirousis and Papadimitriou [KP86]. From [EST94], the search number is informally defined in terms of pebbling to be the minimum number of searchers needed to capture a fugitive who is allowed to move with arbitrary speed about the edges of the graph. For node search number, a searcher blocks all neighboring nodes without the need to move along an incident edge.

Kinnersley in [Kin92] has shown that the pathwidth of a graph is identical to the vertex separation of a graph. The concept of pathwidth has been popularized by the theories of Robertson and Seymour (see for example, [RS85]). Thus, since the gate matrix layout cost, another well-studied VLSI layout problem [KL94, Möh90], equals the pathwidth plus one [FL89], it also equals the vertex separation plus one.

This paper shows that vertex separation is also related to another area besides computer science, namely computational biology.

## 2 Main Result

In this section, we formally define our fixed-parameter problems $k$-ICG and $k$-CVS and then show that they are indeed equivalent.

**Definition 1:** A *layout* $L$ of a graph $G = (V, E)$ is a one to one mapping $L : V \to \{1, 2, \ldots, |V|\}$.

If the order of a graph $G = (V, E)$ is $n$, we conveniently write a layout $L$ as a permutation of the vertices $(v_1, v_2, \ldots, v_n)$. For any layout $L = (v_1, v_2, \ldots, v_n)$ of $G$ let $V_i = \{v_j \mid j \leq i \text{ and } (v_j, v_k) \in E \text{ for some } k > i\}$ for each $1 \leq i \leq n$.

**Definition 2:** The *vertex separation* of a graph $G$ with respect to a layout $L$ is $vs(L, G) = \max_{1 \leq i \leq |G|} \{|V_i|\}$. The *vertex separation* of a graph $G$, denoted by $vs(G)$, is the minimum $vs(L, G)$ over all layouts $L$ of $G$.

The $k$-*coloring* of a graph $G = (V, E)$ is a mapping $color : V \to \{1, 2, \ldots, k\}$. For any subset $V' \subseteq V$, let $Colors(V') = \{color(v) \mid v \in V'\}$.

**Definition 3:** A *colored layout* $L$ of a $k$-colored graph $G = (V, E)$ is layout $L$ such that for all $1 \leq i < n$, $color(v_{i+1}) \notin Colors(V_i)$.

**Problem 4:** <u>COLORED VERTEX SEPARATION (CVS)</u>
*Input:* A $k$-colored graph $G$.
*Parameter:* $k$
*Question:* Is there a colored layout $L$ of $G$ where $vs(L, G) < k$?

**Problem 5:** <u>Intervalizing Colored Graphs (ICG)</u>
*Input:* A $k$-colored graph $G = (V, E)$.
*Parameter:* $k$
*Question:* Is there a properly colored supergraph $G' = (V', E')$ of $G$, $E \subseteq E'$, such that $V = V'$ and $G'$ is an interval graph?

Figure 2 below shows a 3-colored graph with an interval supergraph represented on the left and a colored vertex separation layout given on the right.

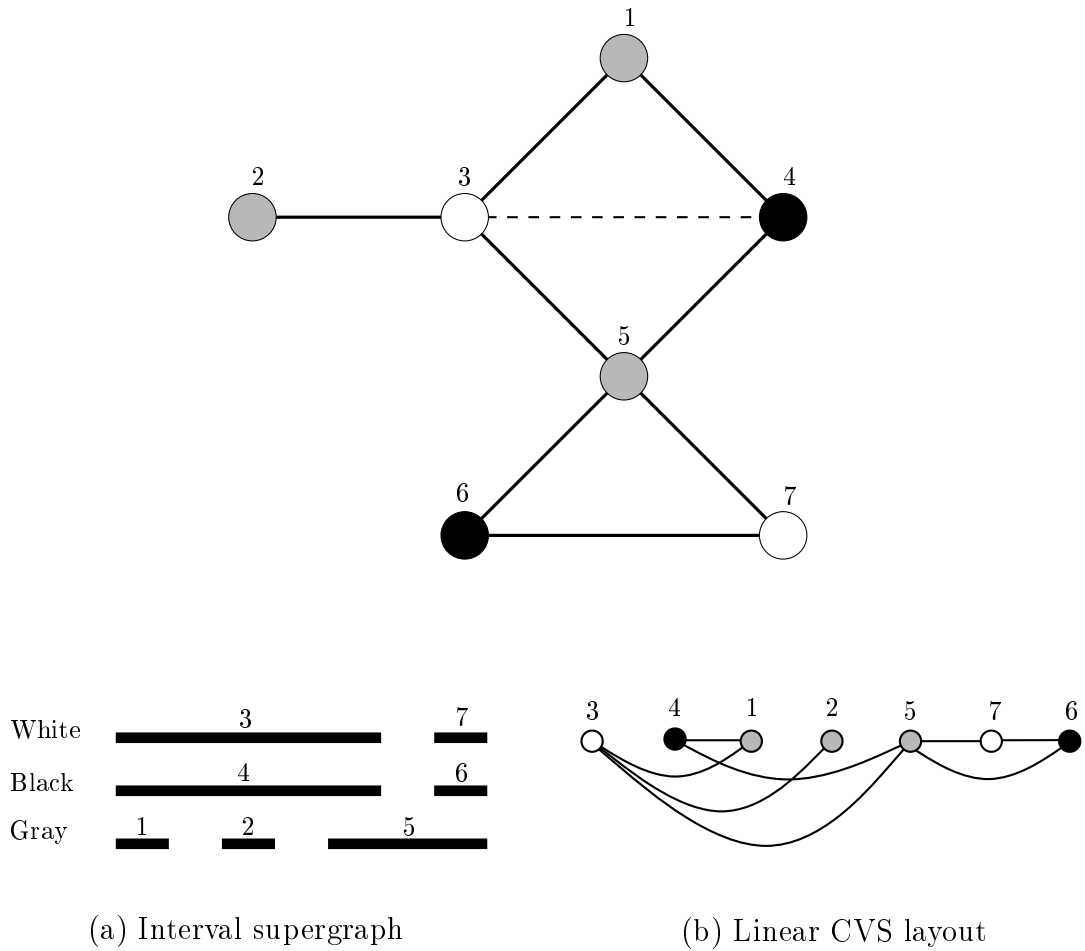

(a) Interval supergraph  (b) Linear CVS layout

Figure 1: Illustrating the $k$-CVS and $k$-ICG problems.

**Theorem 6:** For any fixed positive integer parameter $k$, both $k$-CVS and $k$-ICG are identical problems.

**Proof.** Let $L = (v_1, v_2, \ldots, v_n)$ be a colored layout of a $k$-colored graph $G = (V, E)$. We show how to construct a properly colored supergraph $G'$ that is also an interval graph. For each vertex $v_i \in V$, define the interval:

$$I_{v_i} = [a_{v_i}, b_{v_i}] = [i, \max\{j \mid (v_i, v_j) \in E \vee j = i\} + 0.5]$$

By definition, if edge $(u, v) \in E$ then $I_u \cap I_v \neq \emptyset$. Let $G' = (V, E')$ where $(v_i, v_j) \in E'$ whenever $I_{v_i} \cap I_{v_j} \neq \emptyset$. It suffices to show that $color(v_i) \neq color(v_j)$ for each edge $(v_i, v_j)$ in $E' \setminus E$. Without loss of generality, assume $i < j$ so that $b_{v_i} > a_{v_j}$. Again by the definition of $I_{v_i}$, there exists a vertex $v_k$ such that $j < k$ and $(v_i, v_k) \in E$. This implies that $v_i \in V_{j-1}$. (This also holds for $i = j - 1$.) Now $L$ is a colored layout so $color(v_j) \notin Colors(V_{j-1})$. Thus, $color(v_i) \neq color(v_j)$. Therefore, $G'$ is a properly-colored intervalizable supergraph of $G$.

For any $k$-colored graph $G = (V, E)$ that satisfies ICG, let $\{I_v \mid v \in V\}$ be an interval graph representation of a supergraph $G' = (V, E')$. Let $a_v < b_v$ be the endpoints of the interval $I_v = [a_v, b_v]$ for vertex $v$. Without loss of generality, assume that $a_u = a_v$ implies $u = v$. Let $L = (v_1, v_2, \ldots, v_n)$ be the unique layout such that $i < j$ if and only if $a_{v_i} < a_{v_j}$. We claim that $L$ is a colored layout of $G'$. To prove this claim, we show that $color(v_{i+1}) \notin Colors(V_i)$, $1 \leq i < n$. If there exists a vertex $u \in V_i$ such that $color(u) = color(v_{i+1})$ then by definition of $V_i$ vertex $u$ must be adjacent to a vertex $v_j$ for some $j > i$. Further, $j > i + 1$ since $(u, v_{i+1})$ would not be a properly colored edge. Since $a_u < a_{v_j}$ and $(u, v_j) \in G'$, we must have $b_u > a_{v_j}$ in order to form an overlap. However, $b_u < a_{v_{i+1}} < a_{v_j}$. This is a contradiction to $j > i$. So $u \notin V_i$ if $color(u) = color(v_{i+1})$. Thus $L$ is a colored layout.

Now suppose that for some $r < s$ there exist two vertices $v_r$ and $v_s$ in $V_i$ with the same color. Since $v_r \in V_i$, there exists a vertex $v_j$ with $j > i$ such that $(v_r, v_j) \in E'$. This implies $v_r \in V_{s-1}$. But this implication contradicts the fact $color(v_s) \notin Colors(V_{s-1})$. So $color(v_r) \neq color(v_s)$. Hence any set $V_i \cup \{v_{i+1}\}$ has at most one vertex of each color. Since there are $k$ colors, each $V_i$ must have $k-1$ or fewer vertices. Thus, $vs(L, G) \leq vs(L, G') < k$. □

# 3  Final Comments

Recently, the corresponding general problem of intervalizing a colored graph to an *unit* interval graph has been shown to be NP-hard (and fix-parameter hard for W[1]) by Kaplan and Shamir [KS93] (also see [GGKS93, KST94]). The good news from Kaplan and Shamir's paper is that for each fixed-parameter $k$ (i.e., $k$ colors) this unit interval problem has a polynomial-time algorithm. It is still unknown if a polynomial-time algorithm exists for $k$-ICG, or equivalently $k$-CVS. It is our hope that understanding the original polynomial-time algorithm for the non-colored vertex separation problem may be of some use [EST87].

A related approach for finding a practical $k$-ICG algorithm is based on the easily seen fact that all colored graphs in the $k$-ICG family have pathwidth less than or equal to $k - 1$. The usual polynomial-time algorithms for these types of bounded pathwidth families are constructed as follows: First find a path-decomposition of width $k - 1$ and then use some type of dynamic programming approach on the graph using its decomposition. The tricky part for $k$-ICG is that $k$-ICG is not finite-state (i.e., not representable by linear/tree automaton) for fixed $k$ and hence conventional algorithmic techniques can not be used [FHW93].

However, just because $k$-ICG is not finite-state, we should not avoid altogether the pathwidth structure of the graphs in this family. For small $k$, Bodlaender and Kloks recently developed an algorithm for recognizing and finding path-decompositions of width $k$ in linear time (see [Bod93, BK91, BK93] and [CDF]).